\documentclass[11pt]{amsart}
\usepackage{amssymb}
\usepackage{url}
\usepackage{enumerate}
\usepackage{tikz}
\newtheorem{lemma}{\bf Lemma}[section] 
\newtheorem{theorem}[lemma]{\bf Theorem}
\newtheorem{corollary}[lemma]{\bf Corollary}
\newtheorem{proposition}[lemma]{\bf Proposition}

\newtheorem{definition}[lemma]{\bf Definition}

\DeclareMathOperator{\Amal}{Amal}
\DeclareMathOperator{\Max}{Max}
\DeclareMathOperator{\M}{M}

\begin{document}
\parskip = 0mm
\title[A Question of Gr\"atzer and Lakser from 1971]{A Question of Gr\"atzer and Lakser from the 1971 {\sl Transactions of the American Mathematical Society} \\ \small Forbidden Images of Amalgamation Bases in Finitely-Generated Varieties of Pseudocomplemented Distributive Lattices}
\author{Jonathan David Farley}
\address{Department of Mathematics, Morgan State University, 1700 E. Cold Spring Lane, Baltimore, MD 21251, United States of America, {\tt lattice.theory@gmail.com}}

\author{Dominic van der Zypen}
\address{Swiss Armed Forces Command Support, CH-3003 Bern, Switzerland, {\tt dominic.zypen@gmail.com}}



\begin{abstract}
Gr\"atzer and Lakser asked in the 1971 {\sl Transactions of the American Mathematical Society} if the pseudocomplemented distributive lattices in the amalgamation class of the subvariety generated by ${\bf 2}^n\oplus{\bf 1}$ can be characterized by the property of not having a $*$-homomorphism onto ${\bf 2}^i\oplus{\bf 1}$ for $1<i<n$.

In this article, this question is answered.

If you want to know the answer, you will have to read it (or skip to the last section).
\end{abstract}

\thanks {The first author would like to thank Dr. Sara-Kaja Fischer for introducing him to the topic of amalgamation bases and the author would like to thank Prof. emer. Dr. G\"unter Pilz (then Prof. Dr. Dr. h.c. Pilz) for inviting Dr. Fischer (then Miss Fischer) to Johannes Kepler Universit\"at Linz.}

\maketitle


\def\Qa{\mathbb{Q}_0}
\def\Qb{\mathbb{Q}_1}
\def\Q{\mathbb{Q}}
\def\card{{\rm card}}
\parskip = 2mm
\parindent = 10mm
\def\Part{{\rm Part}}
\def\P{{\mathcal P}}
\def\Eq{{\rm Eq}}
\def\cld{Cl_\tau(\Delta)}
\def\Csing{{\mathcal C}_{\{*\}}}
\def\Cftwo{{\mathcal C}_{{\rm fin}\rangle1}}
\def\Cinf{{\mathcal C}_{\infty}}
\def\Pcf{{\mathcal P}_{\rm cf}}
\def\Fn{{\mathcal F}_n}
\def\proof{{\it Proof. }}
\def\tGeo{{\text{Elena}}}
\def\Geo{{\text{George}}}
\def\tgEO{{\text{eLENA}}}
\def\gEO{{\text{gEORGE}}}
\def\qed{\hfill{}$\Box$}


\vspace*{-4mm} 

\section{Definitions}

For terms and notation not defined, see \cite{DavPriJB}, \cite{GraAJ}, \cite{BalDwiGD} or \cite{McKMcNTayHF}.

A {\it pseudocomplemented distributive lattice} is an algebra $(L,\vee,\wedge,0,1,*)$ such that $(L,\vee,\wedge,0,1)$ is a bounded distributive lattice and $*$ a unary operator---{\it pseudocomplementation}---such that, for all $x,y\in L$, $x\le y^*$ if and only if $x\wedge y=0$. We will call homomorphisms of these algebras $*$-{\it homomorphisms}.

The class of pseudocomplemented distributive lattices is a variety, $\mathcal B_\omega$.  Each subvariety besides $\mathcal B_\omega$ is either $\mathcal B_{-1}$ (the trivial variety) or $\mathcal B_n$ for $n<\omega$, a variety generated by ${\bf 2}^n\oplus{\bf 1}$, an algebra we call $B_n^+$ \cite[Theorems 1 and 6]{LeeGJ}.

Let $\mathcal V$ be a variety, and let $A$ be in $\mathcal V$.  It is an {\it amalgamation base} of $\mathcal V$ if for all $B_0$ and $B_1$ in $\mathcal V$ and for all one-to-one homomorphisms $\alpha_0:A\to B_0$ and $\alpha_1:A\to B_1$, there exist $C$ in $\mathcal V$ and one-to-one homomorphisms $\beta_0:B_0\to C$ and $\beta_1:B_1\to C$ such that $\beta_0\circ\alpha_0=\beta_1\circ\alpha_1$.  The {\it amalgamation class} of $\mathcal V$, $\Amal(\mathcal V)$, is the family of all amalgamation bases of $\mathcal V$.

Clifford Bergman (whose co-discovered theorem on the cancellation of exponents is one of the most beautiful in the theory of ordered sets \cite{BerMcKNagHB}) proved that, for $n<\omega$, every member of $\Amal(\mathcal B_n)$ is a subdirect product of algebras from the set $\{{\bf 2},{\bf 3},{\bf 2}^n\oplus{\bf1}\}$ \cite[Proposition 5.3]{BerHE}, where for $k\in\mathbb N_0$, ${\bf k}$ is the $k$-element chain, $\overline{k}$ the $k$-element antichain, and $\oplus$ denotes the ordinal sum of posets. (For a set $S$, let $\overline S$ be the set $S$ ordered as an antichain.)

Universal algebra pioneer Gr\"atzer (whom Cz\'edli called one of ``the leading experts of lattice theory and, also, of universal algebraíí \cite{FucBeeCzeJH}) and his former student Lakser proved in the 1971 {\sl Transactions of the American Mathematical Society} that if $A$ is in ${\mathcal B}_n$ and has no $*$-homomorphism onto ${\bf 2}^i\oplus{\bf 1}$ for $1<i<n$, then $A$ is in $\Amal({\mathcal B}_n)$ \cite[\S4 Corollary]{GraLakGA}.

They proved that if $A$ is a {\sl finite} pseudocomplemented lattice in ${\mathcal B}_n$ ($n<\omega$), then it is in $\Amal({\mathcal B}_n)$ if {\sl and only if} it has no $*$-homomorphism onto ${\bf 2}^i\oplus{\bf 1}$ when $1<i<n$ \cite[Theorem 6]{GraLakGA}.

In the final paragraph of the final section of their 1971 article in the {\sl Transactions of the American Mathematical Society}, Gr\"atzer and Lakser asked, ``Does Theorem 6 hold for infinite pseudocomplemented lattices?'' \cite[p. 357]{GraLakGA} (Gr\"atzer and Lakser add, ``If not, what is an intrinsic characterization of the algebras in $\Amal({\mathcal B}_n)$ for $2<n<\omega$?'')

We answer this question in $\S5$. 

Gr\"atzer and Lakser stated, ``It seems reasonable, though, that Theorem 6 holds for infinite algebras as well.''  Which way does it go?  The first author went back and forth on this question more times than the top at the end of {\sl Inception}.

We use Priestley duality (\cite{DavPriIC} and the dual of the duality of \cite{PriGE}---prime filters and clopen up-sets versus prime ideals and clopen down-sets).  

An ordered topological space $P$ is {\it totally order-disconnected} if, for all $x,y\in P$ such that $x\nleq y$, there exists a clopen up-set $U$ such that $x\in U$ but $y\notin U$.  A {\it Priestley space} is a compact, totally order-disconnected space.  


Let $\bf P$ be the category of Priestley spaces$+$continuous, order-preserving maps. Let $\bf D$ be the category of bounded distributive lattices$+\{0,1\}$-homomorphisms.  There is a dual equivalence between these two categories \cite{PriGJ}. Our being able to answer the question of Gr\"atzer and Lakser is a testament to the power of Priestley duality. (We are mere altar boys.)

For $L$ in $\bf D$, $P(L)$ denotes its Priestley dual space; for $P$ in $\bf P$, $D(P)$ denotes its Priestley dual bounded distributive lattice. Thus, for $n<\omega$,
$$
P({\bf 2}^n\oplus{\bf 1})\cong {\bf1}\oplus\overline{n}\text{,}
$$
\noindent a Priestley space we call $V_n$. 

We use similar notation for dual morphisms.

If $P$ and $Q$ are objects in $\bf P$ and $f:P\to Q$ a morphism, then $f$ is onto if and only if $D(f)$ is one-to-one; $f$ is an order-embedding if and only if $D(f)$ is onto \cite[Theorem 11.31]{DavPriJB}.  

Hence $L$ will be a subdirect product of a family $(L_{\alpha})_{\alpha\in A}$ of pseudocomplemented distributive lattices with subdirect embedding $e: L\to\prod_{\alpha\in A} L_{\alpha}$ and projections $\pi_{\beta}: \prod_{\alpha\in A} L_{\alpha}\to L_{\beta}$, $(\beta\in A)$ if and only if the duals---$X, (X_{\alpha})_{\alpha\in A}, f:\coprod_{\alpha\in A}X_{\alpha}\to X,$ and $\rho_{\beta}: X_{\beta}\to\coprod_{\alpha\in A} X_{\alpha}$---are such that $f$ is onto and $f\circ\rho_{\beta}:X_{\beta}\rightarrow X$ is an order-embedding for all $\beta\in A$.

But what {\sl are} the Priestley duals of $*$-homomorphism?  For $P$ a poset and $x\in P$, $\Max P$ is the set of maximal elements; let
$$
\M_P(x):=(\Max P)\cap\uparrow x.
$$

For $Y,Y_1$ Priestley duals of pseudocomplemented distributive lattices and $g:Y\to Y_1$ in $\bf P$, then $D(g):D(Y_1)\to D(Y)$ is a $*$-homomorphism if and only if for all $y\in Y$
$$
g[\M_Y(y)]=\M_{Y_1}\big(g(y)\big) \quad\text{\cite[Proposition 3]{PriGE},\cite[p. 47]{DavPriIC}}.
$$

And what are the duals of pseudocomplemented distributive lattices? They are Priestley spaces $P$ such that if $R$ is a clopen up-set of $P$, then $\downarrow R$ is clopen (\cite[Proposition 1]{PriGD}, the dual of \cite[Exercise 11.21]{DavPriJB}).  We will call such spaces $p${\it -spaces} and such maps $p${\it -morphisms}.

By Priestley duality, congruences $\Theta$ of $L=D(Y)$ ($Y$ a $p$-space) correspond to open sets $\Theta'$ of $Y$ such that
$$
\downarrow(\Theta'\cap\Max Y)\subseteq\Theta'
$$
\noindent (\cite{AdaGC}, quoted in \cite[Proposition 13]{PriGE}).   The correspondence is that for all $A_1,A_2\in D(Y)$,
$$
A_1\equiv A_2 \text{ ($\Theta$) if and only if }A_1\setminus\Theta'= A_2\setminus\Theta'
$$
\noindent \cite[\S4, p. 223]{PriGE}.

Using the exponential notation of \cite[p. 67]{FarIFc}, for any poset $P$, $\beta_{\bf 2}(P)$ is the Priestley space (unique up to order-homeomorphism) such that ${\bf2}^P\cong{\bf2}^{\beta_{\bf 2}(P)}$.  We just need to know that if $S$ is a set, then ${\beta_{\bf 2}(\overline S)}$ is a compact, totally disconnected topological space containing $S$ such that $\{s\}$ is open for all $s\in S$ (and we can give the space the antichain ordering).  See \cite[Corollary 3.4, Proposition 3.5, Proposition 3.6, Theorem 3.12, Corollary 3.14, and Proposition 3.15]{FarIFc}.  

\begin{lemma} Let $I$ be a set.  Let $Q$ be a finite poset. Let $\{Q_i\mid i\in I\}$ be a family of pairwise disjoint posets such that for all $i\in I$ $Q_i\cong Q$.  Let $\sum_{i\in I} Q_i$ be the disjoint sum of the posets of the family $\{Q_i\mid i\in I\}$.

Then $\beta_{\bf2}(\sum_{i\in I} Q_i)$ is order-homeomorphic to $\beta_{\bf2}(\overline I)\times Q$.
\end{lemma}

\proof Note that $\beta_{\bf2}(\overline I)\times Q$ is in $\bf P$ and
$$
{\bf 2}^{\sum_{i\in I}Q_i}\cong{\bf 2}^{\overline I\times 
Q}\cong({\bf2}^{\overline I})^Q\cong{\bf2}^{\beta_{\bf2}({\overline I})\times Q}.\quad\text{\cite[Proposition 3.16]{FarIFc}}
$$
\qed

\noindent Also see \cite{CorGF} and \cite{KouSicIA}.

See \cite[Exercise 11.8]{DavPriJB}.

\begin{lemma} Let $Q$ and $R$ be posets.  Then $\beta_{\bf2}(Q+R)\cong\beta_{\bf2}(Q)+\beta_{\bf2}(R)$ (where the second ``+'' denotes the coproduct in $\bf P$).
\end{lemma}

\proof We have
$$
{\bf2}^{Q+R}\cong{\bf2}^Q\times{\bf2}^R\cong{\bf2}^{\beta_{\bf2}Q}\times{\bf2}^{\beta_{\bf2}R}\cong{\bf2}^{\beta_{\bf2}Q+{\beta_{\bf2}R}}. 
$$
\qed

\section{What Is ``Behind the Pulpit''?}
Let $X:=\beta_{\bf2}(\overline{\mathbb N})$.  For $x\in X$, let $\tilde Q_x$ be the poset $\{\tGeo_x\}\oplus\overline{\{\tilde a_x,\tilde b_x,\tilde c_x\}}$ (Figure 2.1).

\begin{center}

    \begin{tikzpicture}[scale=.85]

    \draw[fill] (0,0) circle (.05cm);
    \draw[fill] (-2,2) circle (.05cm);
    \draw[fill] (0,2) circle (.05cm);
    \draw[fill] (2,2) circle (.05cm);

    \draw (0,0) -- (-2,2);
    \draw (0,0) -- (0,2);
    \draw (0,0) -- (2,2);

    \draw (0,-0.5) node {$\tGeo_x$};
    \draw (-2,2.5) node {$\tilde a_x$};
    \draw (0,2.5) node {$\tilde b_x$};
    \draw (2,2.5) node {$\tilde c_x$};

    \draw (0,-2) node {\bf Figure 2.1. The poset $\tilde Q_x$.};

    \end{tikzpicture}
    
\end{center}

Let $\tilde Q:=\sum_{x\in X}\tilde Q_x$, the disjoint sum of the pairwise disjoint posets $\tilde Q_x$ ($x\in X$).  Give $\tilde Q$ the topology with the subbasis for its open sets
$$
\big\{\{\tilde a_u\mid u\in U\},\{\tilde b_u\mid u\in U\},\{\tilde c_u\mid u\in U\},\{\tGeo_u\mid u\in U\}\mid U\text{ open in }X\big\}.
$$
\noindent Note that this makes $\tilde Q$ a Priestley space order-homemorphic to $\overline X\times({\bf1}\oplus\overline 3)$ where $X$ is ordered as an antichain and ${\bf1}\oplus\overline 3$ has the discrete topology.

For $z\in \mathbb N$, let $Q_z$ be $\{\Geo_z\}\oplus\overline{\{a_z,b_z,c_z\}}$ (Figure 2.2).  For $x\in X\setminus{\mathbb N}$, let $Q_x$ be $\{\Geo_x\}\oplus\{a_x,b_x\}$ and let $c_x:=a_x$.

\begin{center}

    \begin{tikzpicture}[scale=.85]

    \draw[fill] (0,0) circle (.05cm);
    \draw[fill] (-2,2) circle (.05cm);
    \draw[fill] (0,2) circle (.05cm);
    \draw[fill] (2,2) circle (.05cm);

    \draw (0,0) -- (-2,2);
    \draw (0,0) -- (0,2);
    \draw (0,0) -- (2,2);

    \draw (0,-0.5) node {$\Geo_z$};
    \draw (-2,2.5) node {$a_z$};
    \draw (0,2.5) node {$b_z$};
    \draw (2,2.5) node {$c_z$};

    \draw[fill] (6,0) circle (.05cm);
    \draw[fill] (4,2) circle (.05cm);
    \draw[fill] (8,2) circle (.05cm);

    \draw (6,0) -- (4,2);
    \draw (6,0) -- (8,2);

    \draw (6,-0.5) node {$\Geo_x$};
    \draw (4,2.5) node {$a_x=c_x$};
    \draw (8,2.5) node {$b_x$};

    \draw (0,-2) node {{\bf Figure 2.2.} The poset $Q_z$ for $z\in\mathbb N$ and $Q_x$ for $x\in X\setminus\mathbb N$.};

    \end{tikzpicture}
    
\end{center}

Let $Q=\sum_{x\in X} Q_x$ (where the summands are pairwise disjoint).  Give $Q$ the topology with subbasis for its {\sl closed} sets
$$
\big\{\{a_u\mid u\in U\},\{b_u\mid u\in U\},\{c_u\mid u\in U\},\{\Geo_u\mid u\in U\}\mid U\text{ clopen in }X\big\}.
$$

Let $\tilde A:=\{\tilde a_x\mid x\in X\}$, $\tilde B:=\{\tilde b_x\mid x\in X\}$, $\tilde C:=\{\tilde c_x\mid x\in X\}$, and $\tgEO:=\{\tGeo_x\mid x\in X\}$, and let
$$
A:=\{a_x\mid x\in X\}\text{, }B:=\{b_x\mid x\in X\}\text{, }C:=\{c_x\mid x\in X\}\text{, and }\gEO:=\{\Geo_x\mid x\in X\}.
$$

Let $\phi:\tilde Q\to Q$ be the map $\tilde a_x\mapsto a_x$, $\tilde b_x\mapsto b_x$, $\tilde c_x\mapsto c_x$, $\tGeo_x\mapsto \Geo_x$ ($x\in X$).

Since $\mathcal {\bf2}^{\overline{\mathbb N}}$ is a Boolean lattice, its Priestley dual space $X$ is totally disconnected, so the family of clopen subets of $X$ is a subbasis {\sl for the closed sets}.

\begin{lemma} Let $\tilde U$ be a clopen subset of $\tilde Q$.  Then $\tilde U\cap\tilde A$ ($\tilde U\cap\tilde B$, $\tilde U\cap\tilde C$, $\tilde U\cap\tgEO$) is clopen in $\tilde Q$ and in the subspace $\tilde A$ (resp., $\tilde B$, $\tilde C$, $\tgEO$), which is homeomorphic to $X$  via the map $x\mapsto \tilde a_x$ (resp., $x\mapsto\tilde b_x$, $x\mapsto\tilde c_x$, $x\mapsto\tGeo_x$) for all $x\in X$.  [In particular, $\tilde A$ (resp., $\tilde B$, $\tilde C$, $\tgEO$) is clopen in $\tilde Q$.]  
\end{lemma}

\proof Since $\tilde A\cup\tilde B\cup\tilde C\cup\tgEO$ is a disjoint union and equals $\tilde Q$ and each set in the union is open in $\tilde Q$, then each set is clopen.  Therefore the intersections in the statement of the lemma are clopen.

Since the aforementioned sets are disjoint, when one considers how open sets are defined from a subbasis, one sees that the subspace topology on $\tilde A$ (resp., $\tilde B$, $\tilde C$, $\tgEO$) is $\big\{\{\tilde a_u\mid u\in U\}\mid U\text{ open in }X\big\}$ (respectively, $\big\{\{\tilde b_u\mid u\in U\}\mid U \text{ open in }X\big\}$, $\big\{\{\tilde c_u\mid u\in U\}\mid U\text{ open in }X\big\}$, and $\big\{\{\tGeo_u\mid u\in U\}\mid U\text{ open in }X\big\}$). \qed

\begin{lemma} The sets $A\cup C$, $B$, and $\gEO$ are clopen in $Q$.  The sets $\{a_z\mid z\in\mathbb N\}$ and $\{c_z\mid z\in\mathbb N\}$ are open in $Q$.

For every clopen subset $U$ of $X$, $\{b_u\mid u\in U\}$ and $\{\Geo_u\mid u\in U\}$ are clopen in $Q$.

For every {\sl open} subset $V$ of $X$, $\{a_v\mid v\in V\setminus\mathbb N\}\cup\{a_v,c_v\mid v\in V\cap\mathbb N\}$ is open in $Q$.

Finally, $Q$ is Hausdorff.
\end{lemma}

\proof Since $Q$ is a disjoint union $(A\cup C)\cup B\cup\gEO$, then $B$ and $\gEO$ are clopen.  Also $A\cup C$ is clopen.

Since $C$ is closed, then $(A\cup C)\setminus C=\{a_z\mid z\in\mathbb N\}$ is open.  Similarly $\{c_z\mid z\in\mathbb N\}$ is open.

Let $U$ be clopen in $X$.  Then $\{b_u\mid u\in U\}$ and $\{b_v\mid v\in X\setminus U\}$ are closed in $Q$ and $B\setminus\{b_v\mid v\in X\setminus U\}=\{b_u\mid u\in U\}$ is open.  The same applies to $\gEO$.

Also, $\{a_v\mid v\in X\setminus U\}$ and $\{c_v\mid v\in X\setminus U\}$ are closed in $Q$.  Hence $(A\cup C)\setminus[\{a_v\mid v\in X\setminus U\}\cup\{c_v\mid v\in X\setminus U\}]$ is open, but it equals
$$
\begin{aligned}
&\{a_v\mid v\in X\setminus\mathbb N\text{ and }v\notin X\setminus U\}\cup\{a_z\mid z\in\mathbb N\text{ and }z\notin X\setminus U\}\cup\{c_z\mid z\in\mathbb N\text{ and }z\notin X\setminus U\}\\
&=\{a_v\mid v\in U\setminus\mathbb N\}\cup\{a_z\mid z\in U\cap\mathbb N\}\cup\{c_z\mid z\in U\cap\mathbb N\}.
\end{aligned}
$$

If $V$ is {\sl open} in $X$, say $V=\bigcup_{\gamma\in\Gamma}U_\gamma$ where $\{U_\gamma\mid \gamma\in\Gamma\}$ is a family of clopen sets $X$.  So
$$
\begin{aligned}
&\{a_v\mid v\in V\setminus\mathbb N\}\cup\{a_z,c_z\mid z\in V\cap\mathbb N\}\\
&=\{a_v\mid v\in(\bigcup_{\gamma\in\Gamma}U_\gamma)\setminus\mathbb N\}\cup\{a_z,c_z\mid z\in(\bigcup_{\gamma\in\Gamma}U_\gamma)\cap\mathbb N\}\\
&=\{a_v\mid v\in(\bigcup_{\gamma\in\Gamma}U_\gamma\setminus\mathbb N)\}\cup\{a_z,c_z\mid z\in\bigcup_{\gamma\in\Gamma}(U_\gamma\cap\mathbb N)\}\\
&=\bigcup_{\gamma\in\Gamma}\{a_v\mid v\in U_\gamma\setminus\mathbb N\}\cup\bigcup_{\gamma\in\Gamma}\{a_z,c_z\mid z\in U_\gamma\cap\mathbb N\}\\
&=\bigcup_{\gamma\in\Gamma}[\{a_v\mid v\in U_\gamma\setminus\mathbb N\}\cup\{a_z,c_z\mid z\in U_\gamma\cap\mathbb N\}]
\end{aligned}
$$
\noindent is open.

Let $q,r\in Q$ be such that $q\ne r$.

If $q,r\in B$, then say $q=b_i$ and $r=b_j$, where $i,j\in X$.  There exists a clopen subset $U$ of $X$ such that $i\in U$ and $j\ne U$.  By the above, $\{b_u\mid u\in U\}$ is a clopen subset of $Q$ containing $q$ but not $r$.

If $q,r\in\gEO$, then say $q=\Geo_i$ and $r=\Geo_j$, where $i,j\in X$.  There exists a clopen subset $U$ of $X$ such that $i\in U$ and $j\ne U$.  By the above, $\{\Geo_u\mid u\in U\}$ is a clopen subset of $Q$ containing $q$ but not $r$.

If $q\in B$ and $r\in Q\setminus B$, then $B$ is a clopen subset of $Q$ containing $q$ but not $r$---similarly if $q\notin B$ but $r\in B$.

If $q\in\gEO$ and $r\notin\gEO$, then $\gEO$ is a clopen subset of $Q$ containing $q$ but not $r$---similarly if $q\notin\gEO$ but $r\in\gEO$.

Now suppose $q,r\in A\cup C$.  If $q,r\in A$---say $q=a_i$ and $r=a_j$ where $i,j\in X$ and $i\ne j$---then there exist disjoint open sets $U,V$ of $X$ such that $i\in U$ and $j\in V$, so
$$
\{a_u\mid u\in U\setminus\mathbb N\}\cup\{a_z,c_z\mid z\in U\cap\mathbb N\}
$$
\noindent and
$$
\{a_v\mid v\in V\setminus\mathbb N\}\cup\{a_z,c_z\mid z\in V\cap\mathbb N\}
$$
\noindent are disjoint open sets of $Q$. The first contains $q$ and the second $r$.

Ditto if $q,r\in C$.

So assume $q\in A\setminus C$ and $r\in C\setminus A$.

We have $q=a_z$ and $r=c_w$ for some $z,w\in\mathbb N$.  Thus $\{a_x\mid x\in\mathbb N\}$ and $\{c_x\mid x\in\mathbb N\}$ are disjoint open sets of $Q$ containing $q$ and $r$, respectively. \qed

\begin{lemma} The subspace $A$ ($B$, $C$, $\gEO$) of $Q$ is homeomorphic to $X$ via the map $x\mapsto a_x$ ($b_x$, $c_x$, $\Geo_x$) for all $x\in X$.
\end{lemma}

\proof The map is a bijection from a compact Hausdorff space to a Hausdorff one (by Lemma 2).  We are done if we can show it is continuous (see, for example, \cite[A.7 Lemma]{DavPriJB}).

Let $U$ be clopen in $X$.  First consider $B$.  By disjointness, the inverse image of $B\cap\{a_u\mid u\in U\}$, $B\cap\{c_u\mid u\in U\}$, or $B\cap\{\Geo_u\mid u\in U\}$ is empty, and the inverse image of $\{b_u\mid u\in U\}$ is $U$.  Thus the inverse image of any closed subbasis member is closed in $X$.  Hence the map is continuous. Ditto for $\gEO$.

Now consider $A$.  Again, let $U$ be a clopen subset of $X$.  The inverse image of $A\cap\{b_u\mid u\in U\}$ and $A\cap\{\Geo_u\mid u\in U\}$ is closed (it is empty) and the inverse image of $\{a_u\mid u\in U\}$ is $U$, which is closed.

Now consider the inverse image of $A\cap\{c_u\mid u\in U\}=\{c_u\mid u\in U\setminus\mathbb N\}$, $U\setminus\mathbb N$.  In $X$, $\mathbb N$ is open, so $U\setminus\mathbb N$ is closed.

Ditto for $C$. \qed

\begin{lemma} The map $\phi$ is continuous and surjective.  Hence $Q$ is compact.
\end{lemma}

\proof Let $U$ be clopen in $X$.  Then $\phi^{-1}(\{a_u\mid u\in U\})$
$$
\begin{aligned}
&=\{\tilde a_u,\tilde c_u\mid u\in U\setminus\mathbb N\}\cup\{\tilde a_u\mid u\in U\cap\mathbb N\}\\
&=(\{\tilde a_u\mid u\in U\})\cup\{\tilde c_u\mid u\in U\}\setminus\{\tilde c_z\mid z\in\mathbb N\}.
\end{aligned}
$$
\noindent Now $\{\tilde a_u\mid u\in U\}$ is clopen in $\tilde A$ by Lemma 1 and $\tilde A$ is clopen in $\tilde Q$ by Lemma 1 so $\{\tilde a_u\mid u\in U\}$ is clopen in $\tilde Q$.  Similarly, $\{\tilde c_u\mid u\in U\}$ is clopen in $\tilde Q$.  As $\{\tilde c_z\mid z\in\mathbb N\}$ is open in $\tilde C$ by Lemma 1 and $\tilde C$ is open $\tilde Q$ by Lemma 1, $\{\tilde c_z\mid z\in\mathbb N\}$ is open in $\tilde Q$.  Hence $\phi^{-1}(\{a_u\mid u\in U\})$ is closed in $\tilde Q$.

Likewise, $\phi^{-1}(\{c_u\mid u\in U\})$ is closed in $\tilde Q$.

Also, $\phi^{-1}(\{b_u\mid u\in U\})=\{\tilde b_u\mid u\in U\}$, which is closed in $\tilde B$ by Lemma 1, and $\tilde B$ is closed in $\tilde Q$ by Lemma 1, so $\{\tilde b_u\mid u\in U\}$ is closed in $\tilde Q$.  Similarly, $\phi^{-1}(\{\Geo_u\mid u\in U\})$ is closed in $\tilde Q$.

Since the sets we took the inverse image of are a subbasis for the closed sets of $Q$, the inverse image of any closed set in $Q$ is closed in $\tilde Q$.

Thus $Q$ is the image of a compact space under a continuous map. \qed

We will use the fact that the image of a closed set of $\tilde Q$ under $\phi$ is a closed set of $Q$ (since $Q$ and $\tilde Q$ are compact and Hausdorff \cite[A.7 Lemma]{DavPriJB}).

\begin{lemma} Let $z\in\mathbb N$.  Then there exists a clopen up-set of $Q$ containing $a_z$ but not $c_z$.  There exists a clopen up-set of $Q$ containing $c_z$ but not $a_z$.
\end{lemma}

\proof Since $\tilde Q$ is a Priestley space, there exists a clopen up-set $\tilde U$ of $\tilde Q$ containing $\tilde a_z$ but not $\tilde c_z$.  Since $\tilde A$ is a clopen up-set of $\tilde Q$ by Lemma 1, we may assume $\tilde U\subseteq \tilde A$ (by taking $\tilde U\cap\tilde A$).

Let $\tilde V:=\{\tilde c_x\mid \tilde a_x\in\tilde U\}\subseteq\tilde C$.  By symmetry, $\tilde V$ is a clopen up-set of $\tilde Q$ containing $\tilde c_z$ but not $\tilde a_z$.

Let $\tilde W:=(\tilde U\cup\tilde V)\setminus\{\tilde c_z\}$.  Since $\{\tilde c_z\}$ is open and $\tilde W\subseteq\Max\tilde Q$, $\tilde W$ is a closed up-set of $\tilde Q$.  Thus $W:=\phi[\tilde W]$ is a closed subset of $Q$, and it is an up-set since $W\subseteq\Max Q$.  It contains $a_z$ but not $c_z$.

Now $K:=\{x\in X\mid \tilde a_x\in\tilde U\}$ is a clopen subset of $X$ by Lemma 1 ($\tilde U$ is a clopen subset of the subspace $\tilde A$).  Hence $X\setminus K$ is a clopen subset of $X$, and hence $\{\tilde a_x\mid x\in X\setminus K\}$ and $\{\tilde c_x\mid x\in X\setminus K\}$ are clopen subsets of the subspaces $\tilde A$ and $\tilde C$, respectively, by Lemma 1.  Since $\tilde A$ and $\tilde C$ are clopen subsets of $\tilde Q$ by Lemma 1, $\{\tilde a_x\mid x\in X\setminus K\}$ and $\{\tilde c_x\mid x\in X\setminus K\}$ are clopen subsets of $\tilde Q$.

Thus $\tilde T:=\{\tilde a_x\mid x\in X\setminus K\}\cup\{\tilde c_x\mid x\in X\setminus K\}\cup\tilde B\cup\tgEO\cup\{\tilde c_z\}$ is a closed subset of $\tilde Q$ (by Lemma 1 and the fact $\tilde Q$ is Hausdorff).  This means $T:=\phi[\tilde T]$ is a closed subset of $Q$, but $T=Q\setminus W$.  Hence $W$ is open.

Therefore $W$ is a clopen up-set containing $a_z$ but not $c_z$. \qed

\begin{lemma} Let $i\in X$.  Then there exists a clopen up-set of $Q$ containing $a_i$ and $c_i$ and not containing $b_i$ and not containing $\Geo_i$.
\end{lemma}

\proof By Lemma 2, $B\cup\gEO$ is a clopen down-set of $Q$ not containing $a_i$ or $c_i$ but containing $b_i$ and $\Geo_i$. \qed

\begin{lemma} Let $i\in X$.  Then there exists a clopen up-set of $Q$ containing $b_i$ but not containing $a_i$, not containing $c_i$, and not containing $\Geo_i$.
\end{lemma}

\proof By Lemma 2, $B$ is a clopen up-set of $Q$. \qed

\begin{lemma} For any clopen subset $U$ of $X$, $\bigcup_{u\in U} Q_u$ is a clopen up-set and down-set of $Q$.
\end{lemma}

\proof Since $\{\tilde a_u\mid u\in U\}$ is a clopen subset of the subspace $\tilde A$ (by Lemma 1) but $\tilde A$ is a clopen subset of $\tilde Q$ (by Lemma 1), $\{\tilde a_u\mid u\in U\}$ is clopen in $\tilde Q$.  Similarly, $\{\tilde b_u\mid u\in U\}$, $\{\tilde c_u\mid u\in U\}$, and $\{\tGeo_u\mid u\in U\}$ are clopen in $\tilde Q$. Hence
$$
\bigcup_{u\in U}\tilde Q_u=\{\tilde a_u\mid u\in U\}\cup\{\tilde b_u\mid u\in U\}\cup\{\tilde c_u\mid u\in U\}\cup\{\tGeo_u\mid u\in U\}
$$
\noindent is clopen in $\tilde Q$.  Thus $\phi[\bigcup_{u\in U}\tilde Q_u]=\bigcup_{u\in U} Q_u$ is closed in $Q$.

Likewise, $\bigcup_{v\in X\setminus U}Q_v$ is closed in $Q$, so $\bigcup_{u\in U} Q_u=Q\setminus\bigcup_{v\in X\setminus U}Q_v$ is open in $Q$. \qed

\begin{lemma} Let $i,j\in X$ be such that $i\ne j$.  Let $q\in Q_i$ and let $r\in Q_j$.  Then there exists a clopen up-set of $Q$ containing $q$ but not $r$.
\end{lemma}

\proof There exists a clopen subset $U$ of $X$ containing $i$ but not $j$.  By Lemma 8, $\bigcup_{u\in U} Q_u$ is a clopen up-set of $Q$ containing $q$ but not $r$. \qed

\begin{lemma} The ordered space $Q$ is a Priestley space.
\end{lemma}

\proof By Lemma 4, $Q$ is compact. 

Let $q,r\in Q$ be such that $q\nleq r$.  Let $q\in Q_i$ and let $r\in Q_j$ where $i,j\in X$.

\quad{\it Case 1. $i\ne j$}

By Lemma 9, there exists a clopen up-set of $Q$ containing $q$ but not $r$.

\quad{\it Case 2. $i=j$}

If $q=a_i$ and $r=b_i$, use Lemma 6.  If $q=a_i$ and $r=c_i$, then $i\in\mathbb N$; use Lemma 5.  If $q=a_i$ and $r=\Geo_i$, use Lemma 6.

If $q=b_i$, use Lemma 7.

If $q=c_i$ and $r=a_i$, then $i\in\mathbb N$; use Lemma 5.  If $q=c_i$ and $r=b_i$, use Lemma 6.  If $q=c_i$ and $r=\Geo_i$, use Lemma 6.

Since $i=j$, we cannot have $q=\Geo_i$.  \qed

\begin{lemma} Let $R$ be a clopen subset of $Q$.  Then $\{x\in X\mid Q_x\cap R\ne\emptyset\}$ is clopen in $X$.
\end{lemma}

\proof Since $A\cap R$ is clopen in the subspace $A$, $\{x\in X\mid a_x\in R\}$ is clopen in $X$.  Similarly, $\{x\in X\mid b_x\in R\}$, $\{x\in X\mid c_x\in R\}$, and $\{x\in X\mid \Geo_x\in R\}$ are clopen in $X$.  

Hence
$$
\begin{aligned}
\{x\in X\mid Q_x\cap R\ne\emptyset\}=\{x\in X\mid a_x\in R\}&\cup\{x\in X\mid b_x\in R\}\cup\{x\in X\mid c_x\in R\}\\
&\cup\{x\in X\mid \Geo_x\in R\}
\end{aligned}
$$

\noindent is clopen in $X$. \qed

\begin{corollary} Let $R$ be a clopen up-set of $Q$.  Then $\downarrow R$ is a clopen subset of $Q$.
\end{corollary}

\proof We see that $\downarrow R=R\cup\{\Geo_x\mid x\in X \text{ and }Q_x\cap R\ne\emptyset\}$.  By Lemma 11 and Lemma 3, $\{\Geo_x\mid x\in X\text{ and }Q_x\cap R\ne\emptyset\}$ is a clopen subset of the subspace $\gEO$.  By Lemma 2, $\gEO$ is a clopen subset of $Q$, so $\{\Geo_x\mid x\in X\text{ and }Q_x\cap R\ne\emptyset\}$ is a clopen subset of $Q$.  Hence $\downarrow R$ is a clopen subset of $Q$. \qed

Consequently, $\phi$ is a $p$-morphism; $\phi$ is surjective and for all $z\in\mathbb N$, $\phi\restriction_{\tilde Q_z}$ is an order-embedding from $\tilde Q_z$ into $Q$ and for all $\tilde q\in\tilde Q$, $\phi[\M_{\tilde Q}(\tilde q)]=\M_Q\big(\phi(\tilde q)\big)$.  Further, as $X\setminus\mathbb N\ne\emptyset$, pick Harry in $X\setminus\mathbb N$.  There is an order-embedding from a poset $V$ order-isomorphic to ${\bf1}\oplus\overline{2}$---say, $\{\text{Fred},\text{Ginger},\text{Ethel}\}$ where $\text{Fred}<\text{Ginger},\text{Ethel}$---onto  $Q_{\text{Harry}}$ such that $\phi[M_V(\text{Fred})]=\phi[\{\text{Ginger},\text{Ethel}\}]=\{a_{\text{Harry}},b_{\text{Harry}}\}=M_Q(\Geo_{\text{Harry}})$ and $\phi[M_V(\text{Ginger})]=\{a_{\text{Harry}}\}$ (say)$=M_Q\big(\phi(\text{Ginger})\big)$ and $\phi[M_V(\text{Ethel})]=\{b_{\text{Harry}}\}=M_Q\big(\phi(\text{Ethel})\big)$.

But Bergmanís theorem about subdirect representations of algebras in $\Amal(\mathcal B_n)$ does not say that every subdirect product of algebras from the set $\{{\bf 2},{\bf 3},{\bf 2}^n\oplus{\bf1}\}$ is in $\Amal(\mathcal B_n)$\dots.

\section{A Necessary Condition: Copying Bergman}

This section really has no business being in the article, as we merely obtain a result that follows from what Bergman already knew (\cite[Theorem 5.4]{BerHE}), which is not surprising, since all we do is copy his proofs.

The following comes from \cite[p. 147]{BerHE}.

\begin{definition} {\rm An algebra $B$ in a variety $\mathcal V$ is {\it congruence extensile} in $\mathcal V$ if, for every algebra $C$ in $\mathcal V$ containing $B$ as a subalgebra, and for every congruence $\Theta$ of $B$, there is a congruence $\Psi$ of $C$ such that $\Psi\cap(B\times B)=\Theta$.}
\end{definition}

We only need a limited version of \cite[Proposition 2.5]{BerHE}. We copy the parts of the proof we need.

\begin{proposition} Let $B$ be an algebra in a variety $\mathcal V$.  Assume that, for every algebra $C$ in $\mathcal V$ containing $B$ as a subalgebra such that $C$ is isomorphic to a product of subdirectly irreducible algebras in $\mathcal V$, and for every congruence $\Theta$ of $B$, there exists a congruence $\Psi$ of $C$ such that $\Psi\cap(B\times B)=\Theta$.

Then $B$ is congruence extensible in $\mathcal V$.
\end{proposition}

\proof Let $D$ be an algebra in $\mathcal V$ containing $B$ as a subalgebra.  Let $C$ be the product of the subdirectly irreducible algebras of $\mathcal V$ in a subdirect product representation of $D$.  View $D$ as a subalgebra of $C$.

Let $\Theta$ be a congruence of $B$.  Then there is a congruence $\Omega$ of $C$ such that $\Omega\cap(B\times B)=\Theta$.  Let $\Psi=\Omega\cap(D\times D)$, a congruence on $D$.  Then $\Psi\cap(B\times B)=\Omega\cap(D\times D)\cap(B\times B)=\Omega\cap(B\times B)=\Theta$. \qed

\begin{lemma} Let $n\ge1$. The only subdirectly irreducibles in $\mathcal B_n$ up to isomorphism are $B_0^+$, \dots, $B_n^+$.
\end{lemma}

\proof By \cite[Theorem 2]{LakGA}, a subdirectly irreducible has the form $B\oplus\bf1$ where $B$ is a Boolean lattice.  If $|B|\ge2^{n+1}$, then we can see from Priestley duality there is an onto $p$-morphism from ${\bf1}\oplus P(B)$ to $V_{n+1}$ \cite[1.32(i) Proposition]{DavPriJB}. 

Hence $B_{n+1}^+$ would be in $\mathcal B_n$ if $B\oplus\bf1$ were, contradicting \cite[Theorems 4 and 5]{LeeGJ} or \cite[Theorem 6]{LeeGJ}. \qed

\begin{lemma} Let $n\ge2$. Let $I$ be a set.  Let $B=\prod_{i\in I} B_n^+$  Assume $C$ is isomorphic to a product of subdirectly irreducible algebras in $\mathcal B_n$ and that $C$ contains $B$ as a subalgebra.

Then each connected component of the underlying poset of $P(C)$ is order-isomorphic to $V_i$ for some $i\in\{0,1,\dots,n\}$.  If $h$ is the Priestley dual of the embedding, then for $0\le i<n$, $h$ sends a component of $P(C)$ order-isomorphic to $V_i$ to a maximal element of $P(B)$.

Each component of the underlying poset of $P(B)$ is the image of a component of the underlying poset of $P(C)$ order-isomorphic to $V_n$ such that $h$ restricted to that component is an order-embedding.
\end{lemma}

\proof The algebra $C$ is isomorphic to $\prod_{i=0}^n \prod_{j\in J_i} B_i^+$, so $P(C)$ is order-homeomorphic to $\sum_{i=0}^n \beta_{\bf2}(\overline{J_i}\times V_i)$ and we know the underlying poset of $\beta_{\bf2}(\overline{J_i}\times V_i)$ is a sum of copies of $V_i$.

If $V$ is one of these copies with $i<n$, then $h\restriction_V$ sends $V$ to a single component $W$ of $P(B)$ order-isomorphic to $V_n$, but $h$ must send the $i$ maximal elements of $V$ to the $n$ maximal elements of $W$ unless the minimum element of $V$ (and hence all of $V$) is sent to a maximal element of $W$.

Since $h$ is the Priestley dual of an embedding, $h$ is onto, so the inverse image of a minimal element of $P(B)$ must come from a component $X\cong V_n$ of $P(C)$---say, $x\in X$ is sent to that minimal element.  Then $x$ cannot be a maximal element, or else $x$ would have to be sent to a maximal element, not a minimal one.  Hence $x$ is a minimal element and the set of maximal elements above it must be sent to the $n$ maximal elements above $h(x)$.  Thus $h\restriction_X$ is an order-embedding. \qed

\begin{lemma} Use the notation of the previous lemma.  Let $\Theta$ be a congruence of $B$.  Then $h^{-1}(\Theta')$ is $\Psi'$ for a congruence $\Psi$ of $C$.
\end{lemma}

\proof Since $h$ is continuous, $h^{-1}(\Theta')$ is open in $P(C)$.  Let $$
m\in h^{-1}(\Theta')\cap\Max P(C).
$$
\noindent  Let $n\in P(C)$ be such that $n\le m$.  Then, by the previous lemma, $h(m)$ is sent to a maximal element of $P(B)$.  Thus $h(m)\in\Theta'\cap\Max P(B)$ and $h(n)\le h(m)$, so, since $\Theta$ is a congruence of $B$, $h(n)\in\Theta'$.  Thus $n\in h^{-1}(\Theta')$. \qed

\begin{proposition} Let $h:X\to Y$ be a $p$-morphism dual to an embedding and let $\Theta$ be a congruence of $D(Y)$.  If $h^{-1}(\Theta')=\Psi'$ for a congruence $\Psi$ of $D(X)$, then for all $B_1,B_2\in D(Y)$,
$$
h^{-1}(B_1)\equiv h^{-1}(B_2) (\Psi) \text{ if and only if }B_1\equiv B_2 (\Theta).
$$
\end{proposition}

\proof We have $B_1\equiv B_2$ ($\Theta$) if and only if $B_1\setminus\Theta'=B_2\setminus\Theta'$ if and only if $h^{-1}(B_1\setminus\Theta')=h^{-1}(B_2\setminus\Theta')$ (we get the backwards direction since $h$ is onto) if and only if $h^{-1}(B_1)\setminus h^{-1}(\Theta')=h^{-1}(B_2)\setminus h^{-1}(\Theta')$ if and only if $h^{-1}(B_1)\equiv h^{-1}(B_2)$ ($\Psi$). \qed

\begin{corollary} Let $n\ge 2$.  Let $I$ be a set.  In the variety $\mathcal B_n$, $\prod_{i\in I} B_n^+$ is congruence extensile.
\end{corollary}

\proof Use Lemma 5, Proposition 6, and Proposition 2. \qed

One finds the following notions in \cite[\S2]{BerHE} and \cite[p. 35]{TayGB}.

\begin{definition} {\rm Let $A$ be a subalgebra of $B$ in a variety $\mathcal V$.  We say $B$ is an {\it essential extension} of $A$ in $\mathcal V$ if, for every non-trivial congruence $\Theta$ of $B$, $\Theta\cap(A\times A)$ is non-trivial.}
\end{definition}

Gr\"atzer and Lakser proved \cite[Lemma 3(a)]{GraLakGA}:

\begin{lemma} An essential extension of a subdirectly irreducible algebra is subdirectly irreducible. \qed
\end{lemma}

We will use part of \cite[Proposition 2.4.2]{BerHE} and \cite[Lemma 3(b)]{GraLakGA}:

\begin{proposition} Let $A$ and $B$ be algebras in a variety $\mathcal V$.  Let $\alpha:A\to B$ be an embedding.

For any congruence $\Theta$ of $B$, there is a congruence $\Psi$ of $B$ such that $\Theta\cap(\alpha[A]\times\alpha[A])=\Psi\cap(\alpha[A]\times\alpha[A])$ and $B$ mod $\Psi$ is an essential extension of the subalgebra $\alpha[A]$ mod $\Psi$. \qed
\end{proposition}

\begin{lemma} Let $n\ge2$.  For $0\le i\le n$, $B_i^+$ embeds in $B_n^+$.
\end{lemma}

\proof Use Priestley duality (or do it directly). \qed

We copy what we need from \cite[Proposition 3.9]{BerHE} and its proof.  We won't even change the notation!  We are not asserting this is new, since we have already mentioned \cite[Theorem 5.4]{BerHE}.

\begin {proposition} Let $n\ge 2$.  Let $A$ be in $\Amal(\mathcal B_n)$.  Then for every $B$ in $\mathcal B_n$ containing $A$ as a subalgebra, every $*$-homomorphism from $A$ to $B_n^+$ can be extended to a $*$-homomorphism from $B$ to $B_n^+$.
\end{proposition}

\proof Let $\alpha:A\to M$ be a $*$-homomorphism, where $M\cong B_n^+$.  Let $\sigma:A\to S$ be a subdirect embedding of $A$ into a product of subdirectly irreducible algebras in $\mathcal B_n$.  By Lemma 11, we actually get an embedding $\beta:A\to C$ into a product of copies of $B_n^+$.  Thus we get embeddings $\iota_A:A\to B$ (the inclusion map) and $\alpha\times\beta:A\to M\times C$.

Since $A$ is in $\Amal(\mathcal B_n)$, there is an algebra $D$ in $\mathcal B_n$ and there are embeddings $\delta:B\to D$ and $\gamma:M\times X\to D$ such that $\delta\circ\iota_A=\gamma\circ(\alpha\times\beta)$.

Let $\Theta$ be the kernel of the projection $\pi$ from $M\times C$ onto $M$.  By Corollary 7, $M\times C$ is congruence extensile in $\mathcal B_n$, so there is a congruence $\Psi$ on $D$ such that $\Psi\cap(\gamma[M\times C]\times\gamma[M\times C])=(\gamma\times\gamma)[\Theta]$---that is, $\{\big((m_1,c_1),(m_2,c_2)\big)\in(M\times C)\times(M\times C)\mid \big(\gamma(m_1,c_1),\gamma(m_2,c_2)\big)\in\Psi\}=\Theta$.

By Proposition 10, $\Psi$ can be chosen so that $D$ mod $\Psi$ is an essential extension of $\gamma[M\times C]$ mod $\Psi$, which is isomorphic to $M$.
By Lemma 9, $D$ mod $\Psi$ is subdirectly irreducible, so, by Lemma 3 and cardinality considerations, $D$ mod $\Psi$ is isomorphic to $B_n^+$ and so equals $\gamma[M\times C]$ mod $\Psi$.  Consider the map
$$
B\to D\to D\text{ mod }\Psi=\gamma[M\times C]\text{ mod }\Psi\to\text{ (via $\gamma$) } M\times C\text{ mod }\Theta\to\text{ (via $\pi$) } M.
$$
\noindent Restricted to $a\in A$, we get
$$
\begin{aligned}
a\mapsto\delta(a)&=\gamma\big((\alpha\times\beta)(a)\big)\\
&=\gamma\bigg(\big(\alpha(a),\beta(a)\big)\bigg)\mapsto\gamma\bigg(\big(\alpha(a),\beta(a)\big)\bigg)\text{ mod }\Psi\\
&\mapsto\big(\alpha(a),\beta(a)\big)\text{ mod }\Theta\mapsto\alpha(a). 
\end{aligned}
$$
\qed

At this point, the reader should be ready use this result to generalize \cite[Theorem 6]{GraLakGA} to infinite lattices. 

\section{Plot Twist: Still Copying Bergman}

We use \cite[Lemma 2]{GraLakGA}; we quote \cite[Lemma 3.3]{BerHE}.

\begin{lemma} Let $A$, $B_0$, and $B_1$ be algebras in a variety $\mathcal V$.  Let $\alpha_0:A\to B_0$ and $\alpha_1:A\to B_1$ be embeddings in $\mathcal V$.

There exist an algebra $C$ and embeddings $\beta_0:B_0\to C$ and $\beta_1:B_1\to C$ such that $\beta_0\circ\alpha_0=\beta_1\circ\alpha_1$ if and only if:

For $\{i,j\}=\{0,1\}$ and $a,b\in B_i$ such that $a\ne b$, there are an algebra $D$ in $\mathcal V$ and homomorphisms $\beta:B_i\to D$, $\gamma:B_j\to D$ such that $\beta(a)\ne\beta(b)$ and $\beta\circ\alpha_i=\gamma\circ\alpha_j$. \qed
\end{lemma}

We copy what we need from Bergman's \cite[Lemma 3.8]{BerHE} and its proof.

\begin{proposition} Let $n\ge2$.  Let $A$ be in $\mathcal B_n$.  Assume that, for every $B$ in $\mathcal B_n$ containing $A$ as a subalgebra, every $*$-homomorphism from $A$ to $B_n^+$ can be extended to $B$.

Then $A$ is in $\Amal(\mathcal B_n)$.
\end{proposition}

\proof Let the situation be as in Lemma 1.  Let $a,b\in B_0$ be distinct.  Writing $B_0$ as a subdirect product of subdirectly irreducible algebras and using the fact that every subdirectly irreducible embeds in $B_n^+$, there is a $*$-homomorphism $\beta:B_0\to D\cong B_n^+$ such that $\beta(a)\ne\beta(b)$.

By the hypothesis, $\beta\circ\alpha_0:A\to D$ can be extended to a $*$-homomorphism $\gamma:B_1\to D$---that is, $\beta\circ\alpha_0=\gamma\circ\alpha_1$. \qed

\begin{lemma} Let $n\ge1$.  Let $Y$ be a $p$-space such that $D(Y)$ is in $\mathcal B_n$.  Then for all $y\in Y$, $|\M(y)|\le n$.
\end{lemma}

\proof This comes from \cite[Theorems 2 and 5]{LeeGJ}. \qed

\begin{theorem} The Priestley dual of the Priestley space $Q$ of \S2 is in $\Amal(\mathcal B_3)$.
\end{theorem}

\proof Let $N=\{w\}\oplus\overline{\{w_1,w_2,w_3\}}\cong V_3$.

Let $Y$ be a $p$-space such that $D(Y)$ is in $\mathcal B_3$ and let $\gamma:Y\to Q$ be a surjective $p$-morphism.

Let $\alpha:N\to Q$ be a $p$-morphism.

Now $\alpha[\{w_1,w_2,w_3\}]\subseteq\Max Q$.

\quad{\it Case 1. $\alpha(w)\in\Max Q$}

In this case, $\alpha[N]=\{\alpha(w)\}\subseteq\Max Q$.

Let $y\in Y$ be such that $\gamma(y)=\alpha(w)$.  Without loss of generality, $y\in\Max Y$.

Define $\beta:N\to Y$ as the constant map with image $\{y\}$.  This is a $p$-morphism.  Then $\gamma\circ\beta=\alpha$.

\quad{\it Case 2. $\alpha(w)=\Geo_z$ for some $z\in\mathbb N$}

Without loss of generality, $\alpha(w_1)=a_z$, $\alpha(w_2)=b_z$, $\alpha(w_3)=c_z$.

Pick $y\in\gamma^{-1}(\Geo_z)$.  Since $|\gamma[\M_Y(y)]|\ge|M_Q\big(\gamma(y)\big)|=|M_Q(\Geo_z)|=3$, we know $y\notin\Max Y$ and $|M_Y(y)|\ge3$.  By Lemma 3, $|M_Y(y)|\le3$, so $M_Y(y)=\{y_1,y_2,y_3\}$.  Without loss of generality, $\gamma(y_1)=a_z, \gamma(y_2)=b_z, \gamma(y_3)=c_z$.

Define $\beta:N\to Y$ by $\beta(w)=y$, $\beta(w_1)=y_1$, $\beta(w_2)=y_2$, $\beta(w_3)=y_3$. This is a $p$-morphism and $\gamma\circ\beta=\alpha$.

\quad{\it Case 3. $\alpha(w)=\Geo_x$ for some $x\in X\setminus\mathbb N$}

Let $Q_x=\{q\}\oplus\overline{\{q_1,q_2\}}$.

Without loss of generality, $\alpha(w_1)=\alpha(w_3)=q_1$, and $\alpha(w_2)=q_2$.

Choose $y\in\gamma^{-1}(q)$.  Since $q_1\ne q_2$, we have $y\notin\Max Y$ and $|\M_Y(y)|\in\{2,3\}$ by Lemma 3.

\qquad{\it Case 3a. $|M_Y(y)|=2$}

Say $M_Y(y)=\{y_1,y_2\}$ where $\gamma(y_1)=q_1$ and $\gamma(y_2)=q_2$.

Define $\beta:N\to Y$ by $\beta(w)=y$, $\beta(w_1)=y_1$, $\beta(w_2)=y_2$, $\beta(w_3)=y_1$.  Then $\beta$ is a $p$-morphism and $\gamma\circ\beta=\alpha$.

\qquad{\it Case 3b. $|M_Y(y)|=3$}

Say $M_Y(y)=\{y_1,y_2,y_3\}$ where $\gamma(y_1)=\gamma(y_3)=q_1$ and $\gamma(y_2)=q_2$.

Define $\beta:N\to Y$ by $\beta(w)=y$, $\beta(w_1)=y_1$, $\beta(w_2)=y_2$, $\beta(w_3)=y_3$.  Then $\beta$ is a $p$-morphism and $\gamma\circ\beta=\alpha$.

By Proposition 2, $D(Q)$ is in $\Amal(\mathcal B_3)$. \qed

\section{Answer to the Question of Gr\"atzer and Lakser from the 1971 {\sl Transactions of the American Mathematical Society}}

No.


\end{document}